\documentclass[a4paper,12pt]{article}
\usepackage{amsthm,amssymb,amsmath}

\title{On the Hausdorff dimension of invariant measures of weakly contracting on average measurable IFS}
\date{}

\theoremstyle{plain}
\newtheorem{thm}{Theorem} 

\newtheorem{lem}[thm]{Lemma}
\newtheorem{cor}[thm]{Corollary}

\newtheorem{example}[thm]{Example}

\numberwithin{equation}{section}

\newcommand{\R}{\mathbb R}
\newcommand{\N}{\mathbb N}
\renewcommand{\epsilon}{\varepsilon}
\renewcommand{\rho}{\varrho}
\renewcommand{\theta}{\vartheta}
\def\dint{\mathop{\displaystyle \int}}
\def\dsum{\mathop{\displaystyle \sum}}

\author{Joanna Jaroszewska\\
\normalsize Institute of Mathematics, Warsaw University \\
\normalsize ul. Banacha 2, 02-097 Warszawa, Poland\\
\normalsize e-mail: asia@venco.com.pl\\
Micha\l\ Rams\\
\normalsize Institute of Mathematics, Polish Academy of Sciences\\
\normalsize ul. \'Sniadeckich 8, 00-950 Warszawa, Poland\\
\normalsize e-mail: rams@impan.gov.pl}
\begin{document}

\maketitle

\begin{abstract}
We consider measures which are invariant under a measurable iterated function system with positive, place-dependent probabilities in a~separable metric space.
We provide an upper bound of the Hausdorff dimension of such a measure if it is ergodic.
We also prove that it is ergodic iff the related skew product is.
\end{abstract}

\def\thefootnote{}

\section{Introduction and statement of a result}

In this note we give a contribution to the study of the multifractal properties of measures which are invariant for iterated function systems.
Recently this aspect of such measures has been widely investigated, e.g. some results concerning their Hausdorff dimension were obtained in \cite{S}, \cite{MS}, \cite{NSB}, \cite{JO}, \cite{R} and \cite{FST}.
For instance, in \cite{S} and \cite{MS}, the systems contracting on average and having Dini-continuous, separated from zero probabilities were considered and the upper bound of the Hausdorff dimension of the unique (in this case) invariant distribution was given.
A~certain class of contracting on average systems with constant probabilities was dealt with in \cite{FST}.
(Note that ''contracting on average'' has different meanings in \cite{FST} and in \cite{S}, \cite{MS}.)
The system of that kind could have more than one probability invariant measure, however, in \cite{FST} the upper estimation of the Hausdorff dimension of any of them was provided.
It seems that it was the first attempt to study iterated function systems without uniqueness of invariant distributions in this respect.
Another such attempt was made in \cite{JO}, where the exact dimension of ergodic invariant measures was calculated for a system which acts on a~compact interval, is non-overlapping and has continuous probabilities.

Here we would like to continue the line of research which we described above.
Assume that $\left(X,\rho\right)$ is a separable metric space and $\{X,S_i,p_i\}, i\in I$ is a finite iterated function system with positive probabilities and with ergodic invariant measure $\mu$.
Given $N\in\N$, denote
\[
h_N(\mu,\delta)= -\dint_{\hspace{-2mm}X} \dsum_{i_1\ldots i_N\in I^N} p_{i_1\ldots i_N}(x) \inf\limits_{y\in B_N(x,i_1\ldots i_N,\delta)} \log p_{i_1\ldots i_N}(y)\,\mu(dx)
\]
where $B_N(x,i_1\ldots i_N,\delta)$ is the set of all points $y\in X$ such that for all $n\in [0,N]$
\[
\rho(S_{i_n}\circ\ldots\circ S_{i_1}(y),S_{i_n}\circ\ldots\circ S_{i_1}(x))<\delta
\]
holds.
We will denote
\[
h_N(\mu)=\lim_{\delta\searrow 0} h_N(\mu, \delta)
\]
and
\[
h(\mu)=\lim_{N\to\infty} \frac 1 N h_N(\mu)
\]

Similarly, let
\[
\lambda_N(\mu,\delta)= \dint_{\hspace{-2mm}X} \dsum_{i_1\ldots i_N\in I^N} p_{i_1\ldots i_N}(x) \sup\limits_{\substack{y\in B_N(x,i_1\ldots i_N,\delta)\\y\neq x}} \log \dfrac {\rho(S_{i_N}\circ\ldots\circ S_{i_1}(x),S_{i_N}\circ\ldots\circ S_{i_1}(y))} {\rho(x,y)}\,\mu(dx)
\]

\[
\lambda_N(\mu)=\lim_{\delta\searrow 0} \lambda_N(\mu, \delta)
\]
and
\[
\lambda(\mu)=\lim_{N\to\infty} \frac 1 N \lambda_N(\mu)
\]
(and we accept $-\infty$ as value of $\lambda(\mu)$). As we will
see, $h_N(\mu,\delta)$ and $\lambda_N(\mu,\delta)$ are monotone
with respect to $\delta$ and subadditive with respect to $N$,
hence for $h(\mu)$ and $\lambda(\mu)$ to exist it is enough to
check that $h_1(\mu,\delta)$ and $\lambda_1(\mu,\delta)$ exist for
some $\delta>0$. Our main result is as follows.

\begin{thm} \label{thm:dim}
Under the assumptions as above,
\begin{equation} \label{eqn:dim}
\dim_H(\mu) \leq  -\frac {h(\mu)} {\lambda(\mu)},
\end{equation}
provided the right-hand side is well defined and nonnegative.
\end{thm}

We would like to emphasize that the assumptions of the above theorem are weaker than the ones that usually appear in the context of place-dependent iterated function systems.
Even the existence of invariant measures is not assured here -- it must be guaranteed by additional assumptions.
However, it seems desirable to strengthen Theorem \ref{thm:dim} so that the field of its applicability would contain the systems with probabilities positive only on certain parts of the space, recently considered by I. Werner (see e.g. \cite{We}).

Note that both $h(\mu)$ and $\lambda(\mu)$ are well known in the case of continuous maps $S_i$ and continuous probabilities $p_i$: $h(\mu)$ is the Kolmogorov-Sinai metric entropy of $\mu$ and $\lambda(\mu)$ is the greatest Lyapunov exponent of the system (with respect to the measure $\mu$).
In such a situation, the formula \eqref{eqn:dim} is a generalization of the well known Hofbauer-Raith formula: ratio of entropy to the Lyapunov exponent.

To satisfy the assumptions of the above theorem, the system must have $\lambda_N(\mu)$ negative or equal $-\infty$ for some $N$ (for otherwise the denominator limit would be from $[0^+,\infty]$ and the whole formula would be negative).
This property, crucial for the proof, could be seen as a weak form of contractibility on average.

Let us present here an example of application of our result.

\begin{example} \label{example}
Let $S_1(x)=x/3$ and $S_2(x)=(x+2)/3$, both maps acting on
$X=[0,1]$. Let $A\subset [0,1]$ be a set with (at most) countable
boundary and let $p\in (0,1/2)$. Set $p_1(x)=p$ for $x\in A$ and
$p_1(x)=1-p$ otherwise. Set $p_2=1-p_1$. This iterated function
system has at least one ergodic invariant measure and every its
ergodic invariant measure $\mu$ satisfies
\begin{equation} \label{eqn:example}
\dim_H(\mu)\leq - \frac {p\log p + (1-p)\log (1-p)} {\log 3}
\end{equation}
\end{example}
The paper is divided as follows.
In the second section we introduce the notation and give introductory information about iterated function systems and Markov operators.
We finish the section with discussion of Example \ref{example}.
We also give there an important result on the relationship between ergodicity of iterated function systems and ergodicity of the corresponding skew product.
This allows us to finish the proof of our main result, which we give in the third section.

\section{Preliminaries}

Let $(X,\rho)$ be a fixed nonempty separable metric space and let $I$ be a finite set of cardinality at least $2$.
The following notation will be used through this paper.
The set of natural numbers will not contain 0, i.e. $\N=\{1,\ldots\}$.
To count elements of covers of $X$, which are needed to estimate the Hausdorff dimension of a~measure, we use the space $\Sigma = I^{\N}$.
We endow it with the product topology of $I$ taken with the discrete metric.
For a sequence $\omega\in\Sigma$, the $n$-th term of $\omega$ is denoted by $\omega_n$, whereas by $\omega^n$ -- the concatenation of the first $n$ terms of $\omega$ (i.e. a finite sequence $(\omega_1,\ldots,\omega_n)\in I^n$).
Such an $\omega^n$ determines the {\it cylinder} $C_{\omega^n}=\{\xi\in\Sigma:\xi^n=\omega^n\}$.
We denote by $\sigma$ the left shift map acting on $\Sigma$, i.e. a map such that $(\sigma \omega)_n=\omega_{n+1}$.
Finally, $|E|$ stands for the diameter of $E\subset X$, whereas $B(x,r)$ denotes the closed ball in $X$ with center at $x$ and radius $r>0$.

Suppose we are given Borel measurable maps $S_i:X\rightarrow X,i\in I,$ and Borel measurable functions $p_i:X\rightarrow [0,1],i\in I,$ such that $\sum_I p_i\equiv 1.$
Then we call the triple $\{X,S_i,p_i\}$ a {\it (measurable) iterated function system}. The functions $p_i$ are called {\it probabilities}.

Iterated function systems are usually studied by means of the corresponding Markov chains.
Generally, if we want to define a discrete-time Markov chain, we can start with fixing a{\it~transition probability function} (t.p.f.) $P:X\times\mathcal{B}(X)\rightarrow [0,1]$, i.e. a~function such that $P(x,\cdot)$ is a probability measure for each $x\in X$ and $P(\cdot,A)$ is a Borel measurable function for each $A\in\mathcal{B}(X)$
(by $\mathcal{B}(Y)$ we denote the family of all Borel subsets of a metric space $Y$).
For example, for a given iterated function system $\{X,S_i,p_i\}$, let us consider the function
\[P:X\times\mathcal{B}(X)\ni (x,A)\mapsto\sum_I p_i(x)1_A(S_i(x)).\]
Clearly it is a~t.p.f. -- we say that it {\it corresponds} to $\{X,S_i,p_i\}$.

Having given a~t.p.f. $P$, we can look at it ''dynamically'': to a~fixed $x\in X$ we assign another point, choosing it randomly -- according to the distribution $P(x,\cdot)$.
To better understand this action, it is convenient to think not about individual points, but about their distributions.
This leads to the definition of a{\it~Markov operator} corresponding to $P$, acting on the set $\mathcal{M}$ of all finite Borel measures on $X$ via the formula
\[\mu P(A)=\int_X P(x,A)\,\mu (dx)\;\textrm{ for }A\in\mathcal{B}(X),\mu\in\mathcal{M}.\]
This operator transforms the set $\mathcal{M}_1=\{\mu\in\mathcal{M}:\|\mu\|=1\}$ of distributions into itself ($\|\cdot\|$ denotes here the total variation norm).

If the action defined above has no influence on a measure $\mu\in\mathcal{M}$, i.e. if $\mu P=\mu$, then $\mu$ is said to be {\it invariant} under $P$.
If additionally $\mu\in\mathcal{M}_1$ is an extremal point (in $\mathcal{M}$) of the set of distributions invariant under $P$ then it is called {\it ergodic}.
Later we use a convenient characterization of ergodic measures in terms of invariant sets -- we call a set $A\in\mathcal{B}(X)$ $\mu$-{\it invariant} provided $P(\cdot,A)=1_A$ $\mu$-a.e., where $\mu$ is a given invariant measure.

Let us now come back to the situation when $P$ corresponds to an iterated function system $\{X,S_i,p_i\}$.
Obviously $\mu P=\sum_I S_{i*}(p_i \mu)$, where $T_*(f\nu)(A)=\int f\cdot 1_A\circ T\,d\nu$ for $\nu\in\mathcal{M},f\in\mathcal{B},A\in\mathcal{B}(X)$ and $T:X\rightarrow X$ is a Borel measurable map.
(By $\mathcal{B}$ we denote the space of all bounded Borel measurable functions on $X$.)
We will say that a~measure is ergodic or a~measure/a~set is invariant under $\{X,S_i,p_i\}$ if $P$ has the appropriate property.

Now we assume that $\{X,S_i,p_i\}$ is the iterated function system with an invariant measure $\mu$.
We are going to construct a measure-preserving transformation which corresponds to the initial system and has similar properties.
For any $x\in X,\omega^n\in I^n$, let
\[p_{\omega^n}(x)=p_{\omega_1}(x)\cdot p_{\omega_2}(S_{\omega_1}(x))\cdot\ldots\cdot p_{\omega_n}(S_{\omega^{n-1}}(x))\]
and
\[S_{\omega^n}=S_{\omega_n}\circ\ldots\circ S_{\omega_1}\]
(we also put $S_{\omega^0}\equiv S_{\omega_0}\equiv {\textrm{id}}_X, p_{\omega_0}\equiv 1$, treating $\omega_0$ as an empty sequence). 
Moreover let $p_x $ be a family of probability measures on $\Sigma$, defined on cylinders in the following way
\[p_x\left(C_{\omega^n}\right)=p_{\omega^n}(x) \textrm{\hspace{3mm} for }\omega^n\in I^n, x\in X. \]
Measures $p_x$ can be, in turn, used to define the probability measure $\nu$ on $X\times \Sigma$ by the formula
\begin{equation}\label{eqn:nu}
\nu(d(x,\omega))=p_x(d\omega) \mu(dx).
\end{equation}
This last measure, $\nu$, is invariant under the skew product $S$ acting on $X\times \Sigma$ as follows
\[S:(x,\omega)\mapsto (S_{\omega_1}(x), \sigma\omega).\]
It is clear that the properties of the aforementioned invariant measures are related.
Later we show equivalence of their ergodicity -- to do this we need some simple though useful facts.

Let us fix a t.p.f. $P$.
\begin{lem} \label{lem:abs}
Assume that $\mu,\mu_1\in\mathcal{M}_1$ are invariant under~$P$, $\mu_1\ll\mu$ and the following condition holds
\begin{equation} \label{eqn:erg}
\mu(A)=0\textrm{\hspace{0.05in} or \hspace{0.05in}}\mu(A)=1\textrm{\hspace{0.05in} for any }A\in\mathcal{B}(X)\;\mu\textrm{-invariant under $P$}.
\end{equation}
Then $\mu_1=\mu$.
\end{lem}

\begin{proof}
Suppose the lemma is false and set $\tau=\mu-\mu_1\neq 0$.
The minimum property of the Jordan decomposition $\tau=\tau^+-\tau^-$ into nonnegative measures $\tau^+$, $\tau^-$ implies that $\tau^+\leq \tau^+ P$ and $\tau^-\leq \tau^- P$.
Hence and because $P$ preserves the total variation norm, the measures $\tau^+,\tau^-$ are invariant under~$P$.
Moreover, according to the Hahn decomposition theorem, there exist disjoint sets $X^+,X^-\in\mathcal{B}(X)$, $X^+\cup X^-=X$, on which the measures $\tau^+,\tau^-$ (respectively) are concentrated.
Now put $X_0=X^+$, $X_m=\{x\in X_{m-1}:P1_{X_{m-1}}(x)=1\}$ for $m\in\N$ and $A=\bigcap_m X_m$.
It is easy to see that $1_A \leq P(\cdot,A)$ -- therefore $A$ is $\mu$-invariant.
Furthermore, $\tau^+(X)=\tau(X^+),\tau^-(X)=-\tau(X^-)$.
As $\tau\neq 0$, both these numbers are positive and so are $\mu(A),\mu(X\setminus A)$ -- the latter is true since $\mu_1\ll\mu$.
This contradicts~\eqref{eqn:erg}.
\end{proof}

\begin{cor} \label{cor:elt}
If $\mu\in\mathcal{M}_1$ is invariant under $P$, then $\mu$ is ergodic iff the condition \eqref{eqn:erg} holds.
\end{cor}

\begin{proof}
Sufficiency of \eqref{eqn:erg} follows from Lemma \ref{lem:abs} whereas necessity is implied by Lemma 1 from \cite{E}.
\end{proof}

Assume now that $\{X,S_i,p_i\}$ is an iterated function system with an invariant distribution $\mu$ and the distribution $\nu$ is defined by \eqref{eqn:nu}.
\begin{lem} \label{lem:erg}
The measure $\mu$ is ergodic w.r.t. $\{S_i,p_i,X\}$ iff the measure $\nu$ is ergodic for $(X\times\Sigma,S)$.
\end{lem}

\begin{proof}
By Corollary \ref{cor:elt} it suffices to prove that $\mu$ satisfies \eqref{eqn:erg} iff $\nu$ is ergodic.
Sufficiency is obvious, so we turn to necessity.
Suppose $\nu$ is not ergodic, i.e. there exists a set $D\in\mathcal{B}(X\times \Sigma)$ which is $S$-invariant and such that
\begin{equation} \label{eqn:pos}
\nu(D)>0,\nu(D^c)>0.
\end{equation}
For any $x\in X,A\in\mathcal{B}(X\times\Sigma)$ let
\begin{gather*}
A_x=\{\omega\in\Sigma:(x,\omega)\in A\},\, \, L_{A}(x)=p_x(A_x),\\
A^X=\{y\in X: L_A(y)>0\}.
\end{gather*}
Clearly $A_x\in\mathcal{B}(\Sigma),A^X\in\mathcal{B}(X)$; the latter is true because $L_A\in\mathcal{B}$.
Indeed, if we set
$\mathcal{L}=\{E\in\mathcal{B}(X\times\Sigma):L_E\in\mathcal{B}\}$,
then $\mathcal{L}$ is a $\lambda$-system containing the $\pi$-system $\mathcal{P}$ of Borel measurable rectangles, subsets of $X\times\Sigma$.
By the Sierpi\'{n}ski-Dynkin theorem on $\pi$-$\lambda$-systems, $\mathcal{L}\supset\sigma(\mathcal{P})=\mathcal{B}(X\times\Sigma)$.

From the S-invariance of $D$ it follows that
\begin{equation} \label{eqn:pld}
\sum_I p_i(x) L_D(S_i(x)) = L_D(x) \textrm{ \, for every }x\in X,
\end{equation}
which, in turn, implies that $P(\cdot, D^X)\leq 1_{D^X}$. Thus $\mu(D^X)=1$; similarly $\mu((D^c)^X)=1$.
Therefore there exists a set $\tilde{X}\in\mathcal{B}(X)$ of full measure $\mu$ such that $(D^c)_x=(D_x)^c \neq \emptyset$ for every $x\in \tilde{X}$.
Consequently, $L_D+L_{D^c}$ is $\mu$-a.e. equal to $1$.
We are going to examine properties of $L_D$ more precisely.

{\it Claim 1. $L_D$ is $\mu$-a.e. constant.}

Set $X_l=L_D^{-1}(l,\infty)$, where $l>0$ is such that $\mu(X_l)>0$.
Consider the measure $\tilde{\mu}$ define by the formula: $\tilde{\mu}(E)=\mu\left(E\cap X_l\right)$ for $E\in\mathcal{B}(X)$.
We are going to prove that $\tilde{\mu}$ is invariant.
Since $\mu$ is so, $\tilde{\mu} P\ll\mu$.
Furthermore, putting
\[g=\frac{d\tilde{\mu}}{d\mu}-\frac{d\tilde{\mu}P}{d\mu}\]
we get $(L_D-l) g\geq 0$ $\mu$-a.e. and $\int g\,d\mu=0$.
Hence
\[\int_X L_D g \, d\mu \geq 0\]
with equality iff $g$ $\mu$-a.e. equal to $0$.
But
\[\int_X L_D g\,d\mu=\int_X L_D\,d\tilde{\mu}-\int_X L_D\,d\tilde{\mu}P=0,\]
where the last equality is a consequence of \eqref{eqn:pld}.
Thus $\tilde{\mu}$ is an invariant measure absolutely continuous w.r.t. $\mu$.
Lemma \ref{lem:abs} gives $\mu=\tilde{\mu}$, which proves our claim.

To finish the proof of Lemma \ref{lem:erg} it suffices to justify

{\it Claim 2. The following disjunction holds}
\begin{equation} \label{eqn:dis}
L_D=1 \hspace{0.2in}\mu\textrm{-a.e.}\hspace{0.4in} \textrm{or} \hspace{0.4in} L_{D^c}=1 \hspace{0.2in}\mu\textrm{-a.e.}
\end{equation}

Assume the contrary and fix an $\epsilon > 0$ so small that the sets $\{L_D>1-\epsilon\}$ and $\{L_{D^c}>1-\epsilon\}$ have both zero measure $\mu$.
Put $\Lambda(A,\omega^n)=\{x\in X:p_x(A_x\cap C_{\omega^n})>(1-\epsilon)p_x(C_{\omega^n})\}$ for any $n\in\N$, $\omega^n\in I^n$ and $A\subset X\times \Sigma$.
We are going to show that
\begin{equation} \label{eqn:den}
\tilde{X} \subset \bigcup_{Z\in\mathcal{F}} Z,
\end{equation}
where $\mathcal{F}=\left\{\Lambda(A,\omega^n):n\in\N,\omega^n\in I^n,A\in\{D,D^c\}\right\}$.

Pick any $x\in \tilde{X}$ then.
Two cases may occur: either $p_x$ is nonatomic or it has at least one atom.
In the second case that atom, let us call it $\omega$, may be a member of $D_x$ (then we conclude that $x\in\Lambda(D,\omega^n)$ for some $n\in\N$) or it may happen that $\omega\in(D^c)_x$ (in this case there exists $n\in\N$ such that $x\in \Lambda(D^c,\omega^n)$).

Now suppose $x\in \tilde{X}$ is such that $p_x$ has no atoms.
By means of the ''Cantor function''-type construction it is easy to build a~measure-preserving homeomorphism between the spaces $(\Sigma,\mathcal{B}(\Sigma),p_x)$ and $([0,1],\mathcal{B}([0,1]),\lambda)$, which transforms cylinders into some intervals.
This and the Lebesgue theorem on density points applied to $D_x$ imply the existence of $\omega^n\in I^n$ such that $x\in\Lambda(D,\omega^n)$.
The proof of \eqref{eqn:den} is finished.

To make use of \eqref{eqn:den} we notice that the family $\mathcal{F}\subset\mathcal{B}(X)$ is countable and each member of $\mathcal{F}$ has zero measure $\mu$.
It is so since, by the definition of $p_x$ and $S$-invariance of sets $D$, $D^c$,
\[\Lambda(A, \omega^n) \subset S_{\omega^n}^{-1} (\{L_A > 1-\epsilon\})\hspace{0.1in}\textrm{ for }A\in\{D,D^c\},n\in\N,\omega^n\in I^n.\]
But this implies the equality $\mu(\tilde{X})=0$ which contradicts the way we chose the set $\tilde{X}$.
\end{proof}

On account of Lemma \ref{lem:erg}, the Egorov theorem and the Birkhoff ergodic theorem, we have an immediate corollary:

\begin{cor} \label{cor:erg}
Let $\mu$ be ergodic.
Then for every $\epsilon>0$ and for any family of Borel measurable functions $\{h_i: X\rightarrow \R\}_I$ satisfying the inequalities
\[-\infty<\int_X \sum_I p_i h_i \, d\mu<0\]
there exist $K>0,A_K\in \mathcal{B}(X\times\Sigma)$ such that $\nu(A_K)<\epsilon$ and
\[\sum_{j=1}^n h_{\omega_j}(S_{\omega^{j-1}}(x)) < K\]
for all $(x,\omega)\in {A_K}^c$ and all $n \in \N$.
\end{cor}

Let us go back to Example \ref{example}.
Let $\mu$ be any invariant measure for $(X, S_i, p_i)$.
As $\mu P^n = \mu$, $\mu([k\cdot 3^{-n}, (k+1)\cdot 3^{-n})) \leq (1-p)^n \searrow 0$, $\mu$ cannot have atoms.
It implies that $\mu(B_\delta(\partial A))\searrow 0$ as $\delta\searrow 0$.

The function integrated in the definition of $h_N(\mu,\delta)$
equals $-N\sum_{i=1}^2 p_i(x) \log p_i(x) \mu(dx)$ for all $x$
whose all trajectories avoid $B_\delta(\partial A)$ for time $N$
and is bounded by $-N\log p$ everywhere. Hence, by Lebesgue
majorized convergence theorem

\[
\frac 1 N h_N(\mu) = -\int_X \sum_{i=1}^2 p_i(x) \log p_i(x) \mu(dx) =- p\log p - (1-p) \log (1-p)
\]
At the same time, $\lambda_N(\mu,\delta)=-N\log 3$ for all $N$ and $\delta$ and \eqref{eqn:example} follows.

The one thing remaining to check is that $(X, S_i, p_i)$ from Example \ref{example} has any invariant measures at all.
Let $\mu_0$ be any probabilistic measure on $X$ and define
\[
\mu_n = \frac 1 n \sum_{m=0}^{n-1} \mu_0 P^m
\]

As $\mu_n$ form a sequence of probabilistic measures on a compact space, they have a subsequence $\mu_{n_k}$ convergent in law to some measure $\mu$.
Let us fix $\epsilon > 0$.
As $\mu_0 P^n([k\cdot 3^{-m}, (k+1)\cdot 3^{-m})) \leq (1-p)^m$ for all $n>m$, the same is true for $\mu$.
It follows that $\mu$ cannot have atoms, hence $\mu(B_r(\partial A)) \leq \frac 3 {13} \epsilon$ for $r$ small enough.
As $\mu_{n_k}$ converges to $\mu$,

\begin{equation} \label{ass3}
\mu_{n_k}(B_r(\partial A) \leq \frac 3 {13} \epsilon
\end{equation}
and

\begin{equation} \label{ass0}
\sum_{D_i} \left|\mu_{n_k}(D_i) - \mu(D_i)\right| \leq \frac 3 {13} \epsilon
\end{equation}
for $k$ big enough (where the sum is taken over the components of the complement of $B_r(\partial A)$).
We will prove that $\mu P = \mu$.

Let $d_{FM}$ be the Fortet-Mourier metric \cite{FM} on the space of finite measures, defined as

\[
d_{FM}(\mu, \nu) = \sup \left|\int_X f d\mu - \int_X f d\nu\right|
\]
where the supremum is taken over Lipschitz functions with Lipschitz constant 1 and with absolute value bounded by 1.
It is well known that the topology defined by Fortet-Mourier metric is equivalent to the usual weak* topology, see \cite{D}.

As $\mu_{n_k}$ converge to $\mu$, we have

\begin{equation} \label{ass1}
d_{FM}(\mu_{n_k}, \mu) \leq \frac 3 {13} \epsilon
\end{equation}
for $k$ big enough.

By the definition of $\mu_{n_k}$ we have

\begin{equation} \label{ass2}
d_{FM}(\mu_{n_k}, \mu_{n_k} P) = \frac 1 {n_k} d_{FM}(\mu_0, \mu_0 P^{n_k}) \leq \frac 1 {n_k} \leq \frac 3 {13} \epsilon
\end{equation}
for $k$ big enough.

We can write

\[
d_{FM}(\mu_{n_k} P, \mu P) \leq  d_{FM}\left((\chi_{B_r(\partial
A)} \mu_{n_k}) P, (\chi_{B_r(\partial A)} \mu) P\right) +
\sum_{D_i} d_{FM} \left((\chi_{D_i} \mu_{n_k}) P, (\chi_{D_i} \mu)
P\right)
\]

As the Fortet-Mourier distance of two measures cannot be greater than the sum of their $L^1$-norms, by \eqref{ass3}

\[
\max (\mu_{n_k}(B_r(\partial A)), \mu(B_r(\partial A))) \leq \frac 3 {13} \epsilon
\]
To estimate the following sum note that on each $D_i$ the iterated function system acts as two linear contracting maps with contraction coefficients $1/3$, one chosen with fixed probability $p$ and the other $(1-p)$.
Hence, the Fortet-Mourier distance of images of two measures is bounded from above by $1/3$ of the Fortet-Mourier distance of the original measures plus the difference of $L^1$-norms of the original measures.
Summing over $D_i$ and applying \eqref{ass0} we get

\[
\sum_{D_i} d_{FM} ((\chi_{D_i} \mu_{n_k}) P, (\chi_{D_i} \mu) P) \leq \frac 1 3 d_{FM}(\mu_{n_k}, \mu) + \sum_{D_i} |\mu_{n_k}(D_i) - \mu(D_i)| \leq \frac 4 {13} \epsilon
\]
Hence,

\[
d_{FM}(\mu_{n_k} P, \mu P) \leq \frac 7 {13} \epsilon
\]
Applying \eqref{ass1} and \eqref{ass2} we get

\[
d_{FM}(\mu, \mu P) \leq \epsilon
\]
As $\epsilon$ was arbitrary, $\mu = \mu P$.

\section{The proof of Theorem \ref{thm:dim}\label{sec:proof}}

First we would like to clear up some simple case. Namely, notice
that w.l.o.g. we may assume that $X$ do not contain any isolated
points. Indeed, every isolated point is either of zero measure
$\mu$ and therefore can be removed from the space (without
changing $\dim_H(\mu)$) or has it positive. In the second case
$\mu$ is concentrated on a~finite set and consequently
$\dim_H(\mu)=0$.

\begin{lem}
If $h_N(\mu,\delta)$ and $\lambda_N(\mu,\delta)$ exist, they are
subadditive in $N$.
\begin{proof}
Both $h_N(\mu,\delta)$ and $\lambda_N(\mu,\delta)$ can be written in the form
\[
\dint_{\hspace{-2mm}X} \dsum_{\omega^N\in I^N} p_{\omega^N}(x) \sup\limits_{y\in B_N(x,\omega^N,\delta)} \phi(x,y,\omega^N)\,\mu(dx),
\]
where
\[
\phi(x,y,\omega^N)=\sum_{n=1}^N \phi(S_{\omega^{n-1}}(x), S_{\omega^{n-1}}(y), \omega_n)
\]
is some real-valued function.
Given $N_1$ and $N_2$, for any fixed $\delta$ we have $y\in B_{N_1+N_2}(x,\omega^{N_1+N_2},\delta)$ if and only if $y\in B_{N_1}(x,\omega^{N_1},\delta)$ and $S_{\omega^{N_1}}(y) \in B_{N_1}(S_{\omega^{N_1}}(x),\omega_{N_1+1}\ldots \omega_{N_1+N_2},\delta)$.
Hence,
\begin{eqnarray*}
\sup\limits_{y\in B_{N_1+N_2}(x,\omega^{N_1+N_2},\delta)} \phi(x,y,\omega^{N_1+N_2}) \leq \sup\limits_{y\in B_{N_1}(x,\omega^{N_1},\delta)} \phi(x,y,\omega^{N_1}) + \\
\sup\limits_{y\in B_{N_2}(S_{\omega^{N_1}}(x),\omega_{N_1+1}\ldots \omega_{N_1+N_2},\delta)} \phi(S_{\omega^{N_1}}(x),y,\omega_{N_1+1}\ldots \omega_{N_1+N_2})
\end{eqnarray*}
As $\mu=\mu P^{N_1}$, $h_{N_1+N_2}(\mu, \delta)\leq h_{N_1}(\mu,
\delta)+h_{N_2}(\mu,\delta)$ for every $\delta>0$ (and analogously
for $\lambda_{N_1+N_2}(\mu, \delta)$).
\end{proof}
\end{lem}

Now let us present the idea of the proof of Theorem \ref{thm:dim}.
We need to prove that $\dim_H(\mu)\leq h_N(\mu, \delta)/\lambda_N(\mu,\delta)$.
We will only give the detailed proof for $N=1$, for higher $N$ one works with $P^N$ instead of $P$ and the proof is almost identical.

We are going to analyze a family of measures $\{\mu_{j,\omega^n}\}$, which sum up to $\mu$ and are associated with a certain finite partition $\{E_j\}$ of $X$.
Every measure $\mu_{j,\omega^n}$ is concentrated on $B(S_{\omega^n}(e_j),\left|S_{\omega^n}(E_j)\right|)$ (where $e_j\in E_j$) -- typically a set of small diameter.
We will choose some of the pairs $(j,\omega^n)$ in such a way that the union of sets $B(S_{\omega^n}(e_j),\left|S_{\omega^n}(E_j)\right|)$ corresponding to the chosen pairs will be both of big measure $\mu$ and geometrically small (see Lemma \ref{lem:lem}).
We will use these balls to estimate the Hausdorff dimension of $\mu$.

We need some additional notations: for $i \in I, x \in X,\theta<0,\delta>0,m \in \N, \omega \in \Sigma$ we will write

\[L^{\delta,\theta}_i(x)= \max\left\{\sup_{\substack{y\in B(x,\delta)\\y\neq x}} \log \frac {\rho(S_i(x),S_i(y))} {\rho(x,y)},\theta\right\}\]

\[H^\delta_i(x)= \inf_{y\in B(x,\delta)} \log p_i(y),\]

\[L^{\delta,\theta}_m(x,\omega) = \frac 1 m \sum_{k=0}^{m-1} L^{\delta,\theta}_{\omega_{k+1}}(S_{\omega^k}(x)),\]

\[H^\delta_m(x,\omega) = \frac 1 m \sum_{k=0}^{m-1} H^\delta_{\omega_{k+1}}(S_{\omega^k}(x)),\]

\[I^{\delta,\theta}_L=\int_X \sum_I p_i L^{\delta,\theta}_i\,d\mu,\]

\[I^\delta_H=\int_X \sum_I p_i H^\delta_i\,d\mu,\]

\[s(\delta,\theta)= I^\delta_H/I^{\delta,\theta}_L.\]

Note that inequality \eqref{eqn:dim} holds if $I^\delta_H=-\infty$ for each $\delta>0$.
So, since $I^\delta_H$ is a nonincreasing function of $\delta$, w.l.o.g. we can assume that there is $\Delta\in(0,1)$ such that for every $\delta\in(0,\Delta)$ the integral $I^\delta_H$ is finite and, at the same time, the corresponding integral from the denominator of the right-hand side of $\eqref{eqn:dim}$ is negative (see comments in the first section) or equal to $-\infty$.
In the latter case there exists a number $\Theta<0$ such that $I^{\Delta,\theta}_L\in(-\infty,0)$ for all $\theta\in(-\infty,\Theta]$.
In the former one we set $\Theta=-\infty$  -- that constant would play no role in the proof then.

Obviously $s$ (considered on $(0,\Delta)\times[-\infty,\Theta]$) is a nondecreasing function of every variable with another one fixed.
Moreover, the limit of $s$, taken as $(\delta,\theta)\rightarrow (0,-\infty)$, is equal to the right-hand side of \eqref{eqn:dim}.

After making these introductory remarks we can start the proof. First we choose $(\delta,\theta)$ from the domain of $s$.
Next we fix $s>s(\delta,\theta)$ and $\iota\in\N$ in an arbitrary way.
We also pick $\epsilon\in (0,1/\iota)$ such that $s>s(\delta,\theta)(1+\epsilon)/(1-\epsilon)$ and then we apply Corollary~\ref{cor:erg} for $\epsilon$ and $\{L^{\delta,\theta}_i\}_I$ -- we are allowed since $I^{\delta,\theta}_L\in(-\infty,0)$.
As a~result we obtain $K,A_K$.

Let $\{E_j\}_J$ be a finite family of nonempty disjoint Borel subsets of $X$, of diameter at most $\delta e^{-K}/4$ and such that
\[\mu\left(\bigcup_J E_j\right)>1-\epsilon,\]
where $J\subset\N$.
We add to this family $E_0=X\setminus \bigcup_J E_j$ to form a partition of $X$.
Then we pick $n\in\N$ big enough ($n>n_0$) to apply Lemma \ref{lem:lem} (see below) and such that
\begin{equation} \label{eqn:n}
\sharp J(\delta e^{-K})^s \exp \left\{n (s(1-\epsilon)I^{\delta,\theta}_L - (1+\epsilon)I^\delta_H)\right\} < 2^{-\iota}.
\end{equation}
Before we formulate this Lemma, we notice that
\[\mu = \sum_{\{J\cup\{0\}\}\times I^m} \mu_{j,\omega^m},\]
where
\[\mu_{j,\omega^m} = S_{\omega^m*} (p_{\omega^m} \mu|_{E_j})\]
and $m\in \N\cup\{0\}$.
If, for any $(j,\omega^m)\in \{J\cup\{0\}\}\times I^m,i\in I$, we put
\[p_i(j,\omega^m)=\frac 1 {\|\mu_{j,\omega^m}\|} \int_X p_i\,d\mu_{j,\omega^m}\]
then
\[\|\mu_{j,\omega^m i}\|= p_i(j,\omega^m) \|\mu_{j,\omega^m}\|;\]
moreover $p_i(j,\omega^m)>0$, as it is an average of $p_i$ under some measure.

\begin{lem} \label{lem:lem}
There exists $n_0\in\N$ such that for every $n>n_0$ from the set $\{J\cup\{0\}\}\times I^n$ one can choose a~subset $Z(\iota)$ satisfying the following conditions:
\begin{itemize}
\item[i)] for $(j,\omega^n)\in Z(\iota)$ and $m\leq n$, $|S_{\omega^m}(E_j)|\leq \delta$,
\item[ii)] $|S_{\omega^n}(E_j)|\leq \delta e^{-K + n (1-\epsilon)I^{\delta,\theta}_L}/2 $ for $(j,\omega^n)\in Z(\iota)$,
\item[iii)] $Z(\iota)$ has at most $\sharp J e^{-n(1+\epsilon) I^{\delta}_H}$ elements,
\item[iv)] $\mu\left(\bigcup_{Z(\iota)} B\!\left(S_{\omega^n}(e_j),\left|S_{\omega^n}(E_j)\right|\right)\right)\geq 1-4\epsilon$ provided $e_j\in E_j$ for all $j\in J$.
\end{itemize}

\begin{proof}
Fix points $e_j\in E_j$.
Since $\mu$ is invariant, we get
\begin{equation} \label{eqn:key}
\mu\left(\bigcup_{Z(\iota)} B\!\left(S_{\omega^n}(e_j),\left|S_{\omega^n}(E_j)\right|\right)\right) \geq 1- \sum_{Z(\iota)^c} \left\|\mu_{j,\omega^n}\right\|,
\end{equation}
regardless of $n$ and $Z(\iota)\subset \{J\cup\{0\}\}\times I^n$.
Hence, we only need to estimate the sum of $\|\mu_{j,\omega^n}\|$ over pairs $(j,\omega^n)$ for which i), ii) or iii) does not hold.

First we are going to establish how big $n_0$ should be.
Lemma \ref{lem:erg} and the Birkhoff theorem show that $L^{\delta,\theta}_m(x,\omega)$ converges to $I^{\delta,\theta}_L$ and $H^{\delta}_m(x,\omega)$ converges to $I^{\delta}_H$ for $\mu$-almost every $x\in X$ and $p_x$-almost every $\omega\in\Sigma$.
By the Egorov theorem we can choose a~big (w.r.t. $\mu$) subset $X_0\in\mathcal{B}(X)$ for which those convergences are uniform for a~big set of $\omega$'s.
Specifying, we have $n_0\in \N,X_0\in\mathcal{B}(X)$ such that
\begin{equation} \label{eqn:x0}
\mu(X_0)\geq 1-\epsilon
\end{equation}
and for every $x\in X_0$, the set of $\omega$ which does not satisfy
\begin{equation} \label{eqn:lh}
|L^{\delta,\theta}_n(x,\omega)-I^{\delta,\theta}_L|+|H^{\delta}_n(x,\omega)-I^{\delta}_H|\leq-\epsilon\max\{I^{\delta,\theta}_L,I^{\delta}_H\}
\end{equation}
for some $n>n_0$, has $p_x$-measure not greater than $\epsilon$.

Now we are going to construct $Z(\iota)$.
We fix $n>n_0$.
Notice that $Z(\iota)$ cannot contain any element of the form $(0,\omega^n)$.
Let us denote the set of all such pairs by $A_1$.
We obtain
\begin{equation} \label{eqn:a1}
\sum_{A_1} \|\mu_{j,\omega^n}\|= \mu(E_0) <\epsilon.
\end{equation}

Next, we want to exclude from $Z(\iota)$ the set $A_2$ of all the pairs $(j,\omega^n) \in A_1^c$ for which i) does not hold.
So we take $(j,\omega^n)\in A_2$, i.e. such that
\[|S_{\omega^m}(E_j)|>\delta\]
for some $m\leq n$, but
\[|S_{\omega^u}(E_j)|\leq\delta\]
for all $u<m$.
We have
\[\delta<\sup_{y\in E_j} \rho(S_{\omega^m}(x),S_{\omega^m}(y))\]
for some $x\in E_j$. Furthermore
\begin{eqnarray*}
\sup_{y\in E_j} \rho(S_{\omega^m}(x),S_{\omega^m}(y)) &\leq &
|E_j| \prod_{u=0}^{m-1}\;\sup_{\substack{y\in S_{\omega^u}^{-1}(S_{\omega^u}\{x\})}^c \cap E_j}
\frac{\rho(S_{\omega^{u+1}}(x),S_{\omega^{u+1}}(y))}
{\rho(S_{\omega^u}(x),S_{\omega^u}(y))}\\
&\leq & |E_j| \exp(mL^{\delta,\theta}_m(x,\omega)),
\end{eqnarray*}
where $\omega$ is any continuation of the finite sequence $\omega^n$.
Thus $E_j\times C_{\omega^n}\subset A_K$, and so
\begin{equation} \label{eqn:a2}
\sum_{A_2} \|\mu_{j,\omega^n}\| <\epsilon.
\end{equation}

The last set to exclude is the one containing the elements for which we cannot make use of the choice of $n_0$.
We distinguish the set $A_3$ of pairs $(j,\omega^n)$ with $j$ such that $E_j$ does not intersect $X_0$;
here, by \eqref{eqn:x0}, we have
\begin{equation} \label{eqn:a3}
\sum_{A_3} \|\mu_{j,\omega^n}\| \leq \epsilon.
\end{equation}
Next, for every $(j,\omega^n)\in A_3^c$ let us choose a~point $x_j\in E_j\cap X_0$.
Denote by $A_4$ the set of pairs $(j,\omega^n)\in A_3^c$ such that \eqref{eqn:lh} does not hold for $(x,\omega)\in \{x_j\}\times C_{\omega^n}$ and $n$.
Again,
\begin{equation} \label{eqn:a4}
\sum_{A_4} \|\mu_{j,\omega^n}\|\leq \epsilon.
\end{equation}

Let $Z(\iota)=\bigcap_{t=1}^4 A_t^c$.
Since i) and \eqref{eqn:lh} hold for $(j,\omega^n)\in Z(\iota)$, we obtain
\begin{equation} \label{eqn:e2}
\prod_{k=0}^{n-1} p_{\omega_{k+1}}(j;\omega^k) \geq
\prod_{k=0}^{n-1} \;\inf_{y\in B(S_{\omega^k}(x_j),|S_{\omega^k}(E_j)|)} p_{\omega_{k+1}}(y) \geq e^{n(1+\epsilon) I^\delta_H}
\end{equation}
and
\[\left|S_{\omega^n}(E_j)\right|\leq
2|E_j|\prod_{k=0}^{n-1} \; \sup_{\substack{y\in B(S_{\omega^k}(x_j),|S_{\omega^k}(E_j)|)\\y\neq S_{\omega^k}(x_j)}} \frac {\rho(S_{\omega^{k+1}}(x_j),S_{\omega_{k+1}}(y))} {\rho(S_{\omega^k}(x_j),y)}\]
which implies that
\begin{equation} \label{eqn:e1}
\left|S_{\omega^n}(E_j)\right|\leq\delta e^{-K + n(1-\epsilon)I^{\delta,\theta}_L}/2.
\end{equation}
By \eqref{eqn:e1}, ii) is satisfied.
Similarly, as $\sum_{\omega^n} \prod_{k=0}^{n-1} p_{\omega_{k+1}}(j;\omega^k)=1$ for every $j\in J$, \eqref{eqn:e2} implies iii).
The assertion iv) follows from \eqref{eqn:a1}, \eqref{eqn:a2}, \eqref{eqn:a3}, \eqref{eqn:a4} and \eqref{eqn:key}.
We are done.
\end{proof}
\end{lem}

The rest of the proof is standard.
Let $e_j\in E_j$ be any points.
The set
\[Y_\iota=\bigcup_{Z(\iota)} B\!\left(S_{\omega^n}(e_j),\left|S_{\omega^n}(E_j)\right|\right)\]
has measure $\mu$ at least $1-4\epsilon$.
At the same time it can be covered with a family
\begin{equation} \label{eqn:cov}
\left\{B\!\left(S_{\omega^n}(e_j),|S_{\omega^n}(E_j)|\right)\right\}_{(j,\omega^n)\in Z(\iota)}
\end{equation}
of at most $\sharp J e^{-n(1+\epsilon)I^{\delta}_H}$ sets of diameter less than $\delta e^{-K+n(1-\epsilon)I^{\delta,\theta}_L}$.

The set
\[Y=\bigcap_{\iota}\bigcup_{\kappa>\iota} Y_{\kappa}\]
has full measure $\mu$.
At the same time $Y$ has zero $s$-dimensional Hausdorff measure, because $Y\subset \bigcup_{\kappa>\iota}Y_\kappa$
for every $\iota$, the diameter of covers \eqref{eqn:cov} converges to $0$ as $\iota\rightarrow\infty$
and \eqref{eqn:n} yields the inequality
\[\sum_{\kappa>\iota} \; \sum_{(j,\omega^n)\in Z(\kappa)} \left|B\!\left(S_{\omega^n}(e_j),\left|S_{\omega^n}(E_j)\right|\right)\right|^s<2^{-\iota}.\]
Since $s>s(\delta,\theta)$ was arbitrary, the Hausdorff dimension of $\mu$ is not greater than $s(\delta,\theta)$.
As $-\theta>0$ can be chosen arbitrarily big, we get the assertion.\qed

\newpage
\bibliography{ref}

\begin{thebibliography}{WWW}

\thispagestyle{empty}

\bibitem[D]{D}
R. M. Dudley,
\newblock {\it Probabilities and Metrics},
\newblock Aarhus Universitet, 1976.

\bibitem[E]{E}
J. H. Elton,
\newblock An ergodic theorem for iterated maps,
\newblock {\it Ergodic Theory Dynam. Systems} 7 (1987), 481--488.

\bibitem[FM]{FM}
R. Fortet, B. Mourier,
\newblock Convergence de la r\'epartition empirique vers la r\'epartition th\'eor\'etique,
\newblock {\it Ann. Sci. École Norm. Sup.} 70 (1953), 267--285.

\bibitem[FST]{FST}
A. H. Fan, K. Simon, H. R. Toth,
\newblock Contracting on average random IFS with repelling fixpoint,
\newblock {\it Journal of Stat. Phys.} 122 (2006), 169--193.

\bibitem[JO]{JO}
A. Johansson, A. \"Oberg,
\newblock Exact dimension of Cantor type measures generated by iterated function systems,
\newblock preprint.

\bibitem[MS]{MS}
J. Myjak, T. Szarek,
\newblock On Hausdorff dimension of invariant measures arising from non-contractive iterated function systems,
\newblock {\it Ann. Mat. Pura Appl.} 181 (2002), 223--237.

\bibitem[NSB]{NSB}
M. Nicol, N. Sidorov, D. Broomhead,
\newblock On the Fine Structure of Stationary Measures in Systems which Contract-on-Average,
\newblock {\it J.~Theor. Probab.} 15 (2002), 715--730.

\bibitem[R]{R}
M. Rams,
\newblock Dimension estimates for invariant measures of con\-tract\-ing-on-average iterated function systems,
\newblock preprint.

\bibitem [S]{S}
T. Szarek,
\newblock The dimension of self-similar measures,
\newblock {\it Bull. Polish Acad. Sci. Math.} 48 (2000), 293--302.

\bibitem[We]{We}
I. Werner,
\newblock Contractive Markov systems,
\newblock {\it J.~London Math. Soc. (2)} 71 (2005), 236--258.

\end{thebibliography}

\end{document}